\documentclass[titlepage,12pt]{article}
\usepackage{latexsym,vatola,amsfonts,amssymb,graphicx,anysize}
\marginsize{25mm}{25mm}{20mm}{25mm}

\begin{document}

\title{
{\LARGE\bf The replacements of signed graphs and Kauffman brackets
 of links\footnote{Project supported by NSFC 10501038}  }
 }
\author{\small Xian'an Jin\footnote{Corresponding author. E-mail: xajin@xmu.edu.cn} \ \ \ Fuji Zhang\\
{\small School of Mathematical Sciences, Xiamen University,}\\
{\small Xiamen, Fujian 361005, P.R.China}\\
}
\maketitle \vskip1cm
\begin{center}
\begin{minipage}{13.5cm}
\begin{center}
{\bf\footnotesize Abstract}
\end {center}
\vskip0.4cm {\ \ \ \footnotesize Let $G$ be a signed graph. Let
$\hat{G}$ be the graph obtained from $G$ by replacing each edge
$e$ by a chain or a sheaf. We first establish a relation between
the $Q$-polynomial of $\hat{G}$[6] and the $W$-polynomial of $G$
[9]. Two special dual cases are derived from the relation, one of
which has been studied in [8]. Based on the one to one
correspondence between signed plane graphs and link diagrams, and
the correspondence between the $Q$-polynomial of signed plane
graph and the Kauffman bracket of link diagram, we can compute the
Kauffman bracket of link diagram corresponding to $\hat{G}$ by
means of the $W$-polynomial of $G$. By this way we use transfer
matrix approach to compute the Kauffman bracket of rational links,
and obtain their closed-form formulae. Finally we provide an
example to point out that the relation we built can be used to
deal with a wide type of link family.
}
\end{minipage}
\end{center}
\vspace*{5mm}
\newpage
\section{\bf\large Introduction}
\noindent

The Kauffman square bracket  $[D]\in\mathbb{Z}[A,B,d]$ of a link
diagram $D$ is defined by the following two rules [1]:
\begin{enumerate}
\item The kauffman square bracket of a diagram consisting of $n$
disjoint simple closed curves in the plane is $d^{n-1}$; \item
$[D]=A[D_v^A]+B[D_v^B]$ for any crossing $v$ of link diagram $D$,
where $D_v^A$ and $D_v^B$ are link diagrams obtained from $D$ by
opening the $A$-channel and $B$-channel respectively.
\end{enumerate}
Kauffman bracket $<D>\in\mathbb{Z}[A]$ of a link diagram $D$ is
related to Kauffman square bracket $[D]$ by
\begin{equation}
<D>=[D]|_{B=A^{-1},d=-A^2-A^{-2}}
\end{equation}

The Kauffman bracket of link diagram is a regular isotopy
invariant of (unoriented) links, the famous Jones polynomial [2]
$V_L(t) \in \mathbb{Z}[t]$ of an oriented link $L$ is related to
the Kauffman bracket by
\begin{equation}
V_L(t)=(-A^3)^{-w(D)}<D>|_{A=t^{-1/4}},
\end{equation}
where $D$ is the diagram of $L$, $w(D)$ and $<$$D$$>$ are the
writhe and the Kauffman bracket of $D$ respectively.

The Kauffman brackets of some link families have been computed,
for example, a recursive formula and a closed-form formula of
Kauffman brackets of pretzel links are obtained in [3] and [4]
respectively.

It is a well-known fact that link diagrams are in one-to-one
correspondence with signed plane graphs. Based on this
correspondence, Kauffman associate the signed graph $G$ a
polynomial $Q[G]\in\mathbb{Z}[A,B,d]$ in three variables $A,B$ and
$d$ [5,6], which specializes to Kauffman square bracket of the
link diagram when $G$ is planar.

Let $D$ be a link diagram corresponding to signed plane graph $G$.
In [7,8], the authors of this paper studied the Kauffman bracket
of the link diagram family \{$D_c$\} corresponding to
homeomorphism class \{$G_c$\} of $G$ obtained from $G$ by
replacing each edge of $G$ by a chain (path) with its sign
reserved, and established a relation between the $Q$-polynomial of
$G_c$ (the Kauffman bracket of $D_c$) and the chain polynomial of
$G$. Dually, the relation between the Kauffman bracket of the link
diagram family \{$G_s$\} corresponding to amallamorphism class
\{$G_s$\} of $G$ obtained from $G$ by replacing each edge of $G$
by a sheaf (multiple edges) with its sign reserved and the sheaf
polynomial of $G$ can also be obtained easily.

In this paper we compute the Kauffman bracket of the link diagram
family \{$\hat{D}$\} corresponding to graphs $\{\hat{G}\}$
obtained from $G$ by replacing some edges by chains and some edges
by sheaves which includes the two special subfamilies \{$D_c$\}
and \{$D_s$\}. We establish a relation between the $Q$-polynomial
of \{$\hat{D}$\} and the Tutte polynomial of $G$ defined by
Bollob$\acute{a}$s [9], which generalize the result in [7,8].

Rational links are a family of links obtained from rational
tangles [10,11] by joining with simple arcs the two upper ends and
the two lower ends. In section 5, we apply the relation to compute
the Kauffman brackets of rational links. By transfer matrix
approach the explicit formula of Kauffman brackets of rational
links are obtained.

The graph $G=(V(G),E(G))$ in this paper allows loops and multiple
edges. We denote by $k(G)$ the number of connected components of
$G$, and by $n(G)=|E(G)|-|V(G)|+k(G)$ the nullity (cyclomatic
number) of $G$. We denote by $E_n$ the graph with $n$ vertices and
no edges. For $F\subset E(G)$ we write $<F>$ for the spanning
subgraph of $G$ with edge set $F$, and $k<F>$ and $n<F>$ for the
number of components and the nullity of this graph respectively.
We use $G-F$ and $G/F$ to denote the graphs obtained from $G$ by
deleting and contracting (that is, deleting the edges and
identifying the ends of each edges in $F$) the edges in $F$
respectively. In particular, for $e\in E$, $G-e$ and $G/e$ denote
the graphs obtained from $G$ by deleting and contracting the edge
$e$ respectively. We refer to [10] for more graph theory.
\section{Link diagram and signed plane graph}
\noindent

We first describe the classical correspondence between link
diagrams and signed plane graphs. A signed graph is a graph with
each edge labelled with a sign ($+$ or $-$). We denote by $s(e)$
the sign of the edge $e$.

The medial graph $M(G)$ of a connected non-trivial plane graph $G$
is a 4-regular plane graph obtained by inserting a vertex on every
edge of $G$, and joining two new vertices by an edge lying in a
face of $G$ if the vertices are on adjacent edges of the face; if
$G$ is trivial (that is, it is an isolated vertex), its medial
graph is a simple closed curve surrounding the vertex (strictly,
it is not a graph); if $G$ is not connected, its medial graph
$M(G)$ is the disjoint union of the medial graphs of its connected
components.

Given a signed plane graph $G$ with medial graph $M(G)$, to turn
$M(G)$ into a link diagram $D=D(G)$, we turn the vertices of
$M(G)$ into crossings by defining a crossing to be over or under
according to the sign of the edge as shown in Figure 1.
Conversely, given a link diagram $D$, shade it checker-boardly
firstly so that the unbounded face is unshaded, then we associate
$D$ with a signed plane graph $G(D)$ as follows: For each shaded
face $F$, take a vertex $V_F$ in $F$, and for each crossing at
which $F_1$ and $F_2$ meet, take an edge $V_{F_1}V_{F_2}$,
furthermore, give each edge $V_{F_1}V_{F_2}$ a sign according to
the type of the crossing as the rules shown in Figure 1. In Figure
1 and Figure 2 below, the broken line is the edge of graph, and
the solid line is the arc of link diagram.
\begin{figure}[htbp]
\centering
\includegraphics[height=1.5cm]{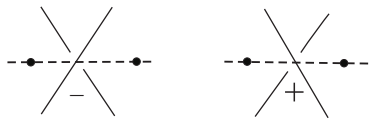}
\caption{The correspondence between a signed edge and a crossing}
\end{figure}

Let $G$ be a signed plane graph. We define the dual signed plane
graph $G^*$ of $G$ as follows: Neglecting the signs of $G$, we
first obtain the dual plane graph of $G$, then give each edge
$e^*$ in the dual plane graph a sign $+$ (resp. $-$) if the sign
of the corresponding edge $e$ in $G$ is $-$ (resp. $+$). Note that
$G$ and $G^*$ correspond to the same link diagram when they are
both connected.

A chain is a graph which is a path and a sheaf is a graph with two
vertices (not necessarily distinct) connected by some parallel
edges. Given a signed plane graph $G$. Let $\hat{G}$ be the graph
obtained from $G$ by replacing each edge $e$ by a chain $c_e$ or a
sheaf $s_e$. By the {\it second Reidemeister Move} shown in Figure
2, the adjacent two edges with different signs in a chain or a
sheaf can cancel each other, without loss of generality we assume
that the signs of the edges in a chain or a sheaf are the same as
that of the edge it replaces. We identify sign $+$ with $+1$, sign
$-$ with $-1$, and define the length of $c_e$ (resp. the width of
$s_e$) as the sum of the signs of the edges in the chain (resp.
the sheaf).

\begin{figure}[htbp]
\centering
\includegraphics[height=5cm]{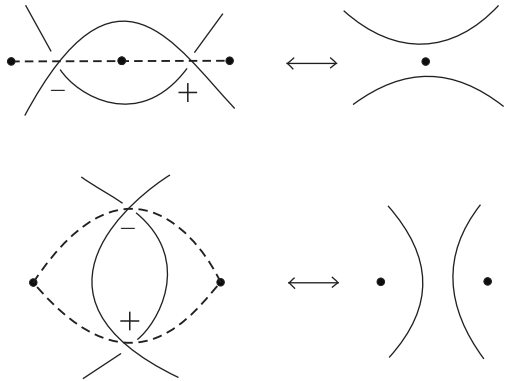}
\caption{The {\it second Reidemeister Move}}
\end{figure}

Let $G$ be a signed plane graph. Since a edge $e$ can be replaced
by a chain or a sheaf and the length of the chain and the width of
the sheaf can be any non-zero integers, we obtain a family of
signed plane graphs $\{\hat{G}\}$ and a link diagram family
\{$\hat{D}$\} accordingly. We call $G$ the associated graph with
the link diagram family \{$\hat{D}$\}. For example, to the graph
$G$ with two parallel edges, we replace one edge by a chain and
the other by a sheaf, then obtain a rational link family (also
called generalized twist link and when $m_2=2$, it is called twist
knot). Figure 3 shows the rational link $m_1m_2$ with $m_1>0$ and
$m_2>0$ and its associated graph. When $m_i<0$, the link diagram
is obtained by replacing the corresponding half-twists by its
mirror image.

\begin{figure}[htbp]
\centering
\includegraphics[height=3.5cm]{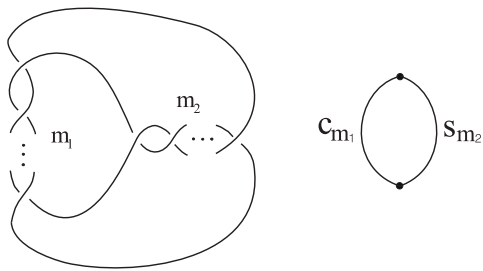}
\caption{Rational link $m_1m_2$ and its associated graph}
\end{figure}

\section{\bf\large Graph polynomials of Tutte type}
\noindent

In this section, we give three graph polynomials of Tutte type
which are all generalizations of Tutte polynomial [13]. They are
$Q$-polynomial for signed graphs, $W$-polynomial for colored
graphs, and chain and sheaf polynomials introduced by Kauffman
[6], Bollob$\acute{a}$s and Riordan[9], and Read and Whitehead
[14] respectively.

\subsection{\bf\large $Q$-polynomial}
\noindent

In [6] Kauffman introduce the Tutte polynomial
$Q[G]=Q[G](A,B,d)\in \mathbb{Z}[A,B,d]$ for signed graph $G$.
Hereafter we assume that $X=A+Bd$ and $Y=Ad+B$. We can redefine
the $Q$-polynomial by the following recursive rules:
\begin{enumerate}
\item
    \begin{equation}
     Q[E_n]=d^{n-1}.
    \end{equation}
\item
   \begin{enumerate}
     \item If $e$ is a loop, then
      \begin{equation}
      \aligned
      &Q[G]=XQ[G-e]\ \ {\rm if}\  s(e)=-,\\
      &Q[G]=YQ[G-e]\ \ {\rm if}\  s(e)=+.
      \endaligned
      \end{equation}
     \item If $e$ is not a loop, then
      \begin{equation}
      \aligned
      &Q[G]=AQ[G-e]+BQ[G/e]\ \ {\rm if}\  s(e)=-,\\
      &Q[G]=BQ[G-e]+AQ[G/e]\ \ {\rm if}\  s(e)=+.
      \endaligned
      \end{equation}
    \end{enumerate}
\end{enumerate}

{\bf Theorem 1}[6] Let $G$ be a signed plane graph. Let $D(G)$ be
the corresponding link diagram. Then $Q[G]=[D(G)]$.

Theorem 1 can be obtained by comparing the rules defining Kauffman
square bracket of link diagram with rules (3)-(5) defining
$Q$-polynomial of the corresponding signed plane graph. On the
other hand, since $G$ and $G^*$ correspond to the same link
diagram, Theorem 1 suggests that $Q[G]=Q[G^*]$.

{\bf Note1:} When $B=A^{-1}$ and $d=-A^2-A^{-2}$, we have
$X=-A^{-3}$ and $Y=-A^3$.

\subsection{\bf\large $W$-polynomial}
\noindent

A colored graph is a graph $G$ together with a function $c$ from
$E(G)$ to an arbitrary set $\Lambda$ of colors, which is a
generalization of signed graph. In [9] Bollob$\acute{a}$s and
Riordan introduced a Tutte polynomial for colored graphs. We
consider a variant for our purpose, and denoted it by $W(G)$. The
polynomial $W(G)=W(G)(t,z_1,z_2)$ for a colored graph $G$ is
defined by the following recursive rules:
\begin{enumerate}
\item \begin{equation}
W(E_n)=t^{n-1}.
\end{equation}
\item Let $c(e)=\lambda$. Then
\begin{enumerate}
\item If $e$ is a bridge, then
\begin{equation}
W(G)=(x_{\lambda}+z_1y_{\lambda}) W(G/e).
\end{equation}
\item If $e$ is a loop, then
\begin{equation}
W(G)=(x_{\lambda}z_2+y_{\lambda})W(G-e).
\end{equation}
\item
 If $e$ is neither a bride nor a loop, then
\begin{equation}
W(G)=x_{\lambda} F(G/e)+y_{\lambda} W(G-e).
\end{equation}
\end{enumerate}
\end{enumerate}

$W$-polynomial is a naturally generalization of $Q$-polynomial
with $\Lambda=\{+,-\}$. The well-definedness of $W(G)$ is based on
the fact that it can be expressed as the sum

{\bf Theorem 2} [9]
$$
W(G)=t^{k(G) -1}\sum_{S\subset E(G)}\{\prod_{e\in
S}x_{c(e)}\}\{\prod_{e\notin
S}y_{c(e)}\}z_1^{k<S>-k(G)}z_2^{n<S>}.
$$

We denote by $G_1\cup G_2$ the union of two disjoint graphs $G_1$
and $G_2$, and by $G_1\cdot G_2$ the union of two graphs with only
one vertex in common. By Theorem 2, we have

{\bf Theorem 3}
\begin{enumerate}
\item $W(G_1\cup G_2)=tW(G_1)W(G_2)$;
\item $W(G_1\cdot
G_2)=W(G_1)W(G_2)$.
\end{enumerate}

A graph $F=(V',E')$ is spanning forest of $G=(V,E)$ if $V'=V,
E'\subseteq E$, and each component of $F$ is a spanning tree of a
component of $G$. Let $G$ be a colored graph and let us consider
an order on its edge set. Let $F$ be a spanning forest of $G$.
Call an edge $e\in E(F)$ internally active if it is the smallest
edge in the unique cut induced by deleting $e$ in $F$ with respect
to the given order and internally inactive otherwise. Similarly,
call an edge $e\notin E(F)$ externally active if it is the
smallest edge in the unique cycle induced by adding $e$ in $F$
with respect to the given order and externally inactive otherwise.
$W(G)$ also has the spanning forest expansion as the various
(weighted) graph polynomials of Tutte type [6,13,15].

{\bf Theorem 4} [9] For any given order on the edge set of $G$, we
have
$$
W(G)=t^{k(G)
-1}\sum_F\{\prod_{IA}(x_{c(e)}+z_1y_{c(e)})\}\{\prod_{EA}(x_{c(e)}z_2+y_{c(e)})\}\{\prod_{II}x_{c(e)}\}\{\prod_{EI}y_{c(e)}\},
$$
where the sum is over all spanning forests of $G$, and the
products are over the edges of internal active, external active,
internal inactive, external inactive respect to $F$ respectively.

\subsection{\bf\large Chain and sheaf polynomials}
\noindent

Chain and sheaf polynomials are introduced by Read and Whitehead
in [14] for studying the chromatic polynomial of homeomorphism
class and the flow polynomial of amallamorphism class of graphs.
They are both the weighted versions of Tutte polynomial [16].
These two dual polynomials are both defined on the graph whose
edges have been labelled with $a,b,c\ldots$ \footnote {We usually
use the label of an edge to represent this edge.}. We denoted the
two polynomials by $Ch[G]$ and $Sh[G]$ respectively.

Before defining the chain and sheaf polynomials, we first give the
definitions of the flow and tension polynomials [13]. Let $G$ be a
graph and let $\overrightarrow{G}$ be an orientation of $G$
obtained by assigning an arbitrary but fixed orientation
$\overrightarrow{e}$ to each edge $e$ of $G$. Let
$\overrightarrow{E}$ denote the set of oriented edges of
$\overrightarrow{G}$. A map $f:\overrightarrow{E}\longrightarrow
\mathbb{Z}_q$ is called a $q$-flow if, for each vertex $v\in
V(G)$, $\sum_+f(\overrightarrow{e})=\sum_-f(\overrightarrow{e})$,
where $\sum_+$ denotes summation over all edges \overrightarrow{e}
incident with $v$ and oriented towards $v$, and $\sum_-$ denotes
summation over all edges \overrightarrow{e} incident with $v$ and
oriented away from $v$. A map $f:\overrightarrow{E}\longrightarrow
\mathbb{Z}_q$ is called a $q$-tension if there is a map
$g:V\longrightarrow \mathbb{Z}_q$ called potential function such
that, for each edge $\overrightarrow{e}=(u,v)\in
\overrightarrow{E}$ oriented from $u$ to $v$,
$f(\overrightarrow{e})=g(u)-g(v)$. $f$ is called nowhere zero if,
for each $\overrightarrow{e}\in \overrightarrow{E}$,
$f(\overrightarrow{e})\neq 0$.

The flow polynomial $F[G](q)$ gives, for each positive integer
$q$, the number of nowhere-zero $q$-flows on $\overrightarrow{G}$.
The tension polynomial $T[G](q)$ gives, for each positive integer
$q$, the number of nowhere-zero $q$-tensions on
$\overrightarrow{G}$. It is not difficult to see that the flow and
tension polynomials do not depend on the orientation of $G$.

The chain and sheaf polynomials can be defined as follows:
$$
Ch[G]=\sum_{Y} {F[G-Y](1-w) \epsilon (Y)}
$$
and
$$
Sh[G]=\sum_{Y} {T[G/Y](1-w) \epsilon (Y)},
$$
where the summations are both over all subsets of $E(G)$; $F[G-Y]$
and $T[G/Y]$ are the flow polynomial of $G-Y$ and the tension
polynomial of $G/Y$; $\epsilon (Y)$ is the product of labels of
the edges in $Y$.

We substitute $1-w$ for the variable of the flow and tension
polynomials in the above definition in order that our definition
is consistent with the original definition by Read and Whitehead.
The following theorem is implicit in [14] and explicit in [16] and
can also be proved directly by dividing the summand in the
definition into two parts, one is over all the $Y$'s which contain
the edge $a$, the other is over all the the $Y$'s which do not
contain the edge $a$.

{\bf Theorem 5}[14,16] The chain polynomial satisfies the
following recursive rules:
\begin{enumerate}
\item
\begin{equation}
Ch[E_n]=1.
\end{equation}
\item
If $a$ is a loop of $G$, then
\begin{equation}
Ch[G]=(a-w)Ch[G-a],
\end{equation}
If $a$ is not a loop, then
\begin{equation}
Ch[G]=(a-1)Ch[G-a]+Ch[G/a].
\end{equation}
\end{enumerate}

Dually, the sheaf polynomial satisfies the following recursive
rules:
\begin{enumerate}
\item
\begin{equation}
Sh[E_n]=1.
\end{equation}
\item
If $a$ is a bridge of $G$, then
\begin{equation}
Sh[G]=(a-w)Sh[G/a]
\end{equation}
If $a$ is not a bridge, then
\begin{equation}
Sh[G]=(a-1)Sh[G/a]+Sh[G-a].
\end{equation}
\end{enumerate}

\section{\bf\large {The relation between} $Q[\hat{G}]$ \bf\large {and} $W(G)$}
\noindent

Let $G$ be a signed graph. Let $\hat{G}$ be the graph obtained
from $G$ by replacing each edge $e$ by a chain $c_e$ or a sheaf
$s_e$. In this section we establish a relation between the
$Q$-polynomial of $\hat{G}$ (hence Kauffman bracket of
$\hat{D}=D(\hat{G})$) and the $W$-polynomial of $G$. Two special
cases are also derived.

{\bf Lemma 6} [8] If $e\in E(G)$ is replaced by a chain with
length $n$, then
\begin{equation}
Q[\hat{G}]={X^n-A^n\over
d}\,Q[\widehat{G-e}]+A^{n}\delta\,Q[\widehat{G/e}],
\end{equation}
where $\delta=d$ if $a$ is a loop of $G$, and $\delta=1$
otherwise.

{\bf Lemma 7} If $e\in E(G)$ is replaced by a sheaf with width
$n$, then
\begin{equation}
Q[\hat{G}]=B^{n}\,Q[\widehat{G-e}]+{Y^{n}-B^{n}\over
d}\delta\,Q[\widehat{G/e}],
\end{equation}
where $\delta=d$ if $e$ is a loop of $G$, and $\delta=1$
otherwise.

{\bf Proof.} We suppose $s(e)=+$ firstly.

{\bf Case 1.} If $e$ is a loop of $G$, by Equation (4) then
$$
\aligned
Q[\hat{G}]&=Y^{n}Q[\widehat{G-e}]\\
 &=B^{n}Q[\widehat{G-e}]+{Y^{n}-B^{n}\over d}d\,Q[\widehat{G/e}].\\
\endaligned
$$

{\bf Case 2.} If $e$ is not a loop of $G$.\\
Suppose that the edges in $s_e$ are $e_1,e_2,\cdots,e_{n}$
successively, by Equations (4) and (5) then
$$
\aligned
Q[\hat{G}]&=BQ[\hat{G}-e_1]+AQ[\hat{G}/e_1]\\
 &=BQ[\hat{G}-e_1]+AY^{n-1}Q[\widehat{G/e}]\\
 &=B^2Q[\hat{G}-e_1-e_2]+(AY^{n-1}+BAY^{n-2})Q[\widehat{G/e}]\\
&=\cdots \cdots\\
&=B^{n_a}Q[\hat{G}-e_1-e_2-\cdots -e_n]+(AY^{n-1}+BAY^{n-2}+\cdots +B^{n-2}AY+B^{n-1}A)Q[\widehat{G/e}]\\
&=B^{n}Q[\widehat{G-e}]+{Y^{n}-B^{n}\over d}\,Q[\widehat{G/e}].
\endaligned
$$

The argument, if $s(e)=-$ is virtually identical to that above,
with $A$ and $B$, $X$ and $Y$, $n$ and $-n$ interchanged
simultaneously. Notice that $B=A^{-1},Y=X^{-1}$, the lemma also
holds. $\Box$
\\

Let $c_n$ be the color of the edge $e\in E(G)$ representing the
corresponding chain with length $n$ in $\hat{G}$ and $s_n$ the
color of the edge $e\in E(G)$ representing the corresponding sheaf
with width $n$ in $\hat{G}$.

{\bf Theorem 8} Let $\hat{G}$ and $G$ be signed graphs explained
above. In $W(G)$, if we set
\begin{equation}
x_{c_n}=A^n, y_{c_n}={X^n-A^n\over d}
\end{equation} and
\begin{equation}
x_{s_n}={Y^n-B^n\over d},y_{s_n}=B^n,
\end{equation} then
\begin{equation}
Q[\hat{G}]=W(G)(d,d,d).
\end{equation}

{\bf Proof:} By Equations (3) and (6), the theorem holds when
$G=E_n$. If $z_1=t$, by Theorem 3, Equation (7) can be replaced by
Equation (9). Comparing the coefficients (16) and (17) with (8)
and (9), the theorem is proved.$\Box$

Now we study two special cases of Theorem 8. Let $G$ be a signed
graph. Let $G_c$ be the graph obtained from $G$ by replacing each
edge $e$ by a chain $c_e$ and $G_s$ be the graph obtained from $G$
by replacing each edge $e$ by a sheaf $s_e$.

By Theorem 8, we have
\begin{equation}
Q[G_c]=W(G)(d,d,d),
\end{equation}
where
$$x_{c_n}=A^n,\ \ \ \ y_{c_n}={X^n-A^n\over d}.$$
Let
$$\tilde{W}(G)=d^{1-p(G)}W(G)(d,d,d),$$
where $p(G)$ is the number of vertices of $G$. We have
\begin{enumerate}
\item
\begin{equation}
\tilde{W}(E_n)=1.
\end{equation}
\item
If $e$ is a loop of $G$, then
\begin{equation}
\tilde{W}(G)={A^n\over d}((X/A)^n-(1-d^2))\tilde{W}(G-e).
\end{equation}
If $e$ is not a loop, then
\begin{equation}
\tilde{W}(G)={A^n\over
d}(((X/A)^n-1)\tilde{W}(G-e)+\tilde{W}(G/e)).
\end{equation}
\end{enumerate}
Comparing Equations (22)-(24) with Equations (10)-(12), we obtain

{\bf Corollary 9}[8] In $Ch[G]$, if we replace $w$ by $1-d^2$, and
replace $a$ by $(X/A)^{n_a}$ for every chain $a$, where $n_a$ is
the length of the chain $a$, then we have
\begin{equation}
Q[G_c]={A^{\sum_{e\in E(G)}n_a}\over d^{q-p+1}}Ch[G],
\end{equation}
where $p$ and $q$ are the number of vertices and edges of graph
$G$, respectively.

Similarly, by Theorem 8, we have
\begin{equation}
Q[G_s]=W(G)(d,d,d)
\end{equation}
where
$$x_{s_n}={Y^n-B^n\over d},\ \ \ \ y_{s_n}=B^n.$$
Let
$$\tilde{W}(G)=d^{1-2k(G)+p(G)}W(G)(d,d,d).$$
We have
\begin{enumerate}
\item
\begin{equation}
\tilde{W}(E_n)=1.
\end{equation}
\item
If $e$ is a bridge of $G$, then
\begin{equation}
\tilde{W}(G)={B^n}((Y/B)^n-(1-d^2))\tilde{W}(G/e).
\end{equation}
If $e$ is not a bridge, then
\begin{equation}
\tilde{W}(G)={B^n}(((Y/B)^n-1)\tilde{W}(G-e)+\tilde{W}(G/e)).
\end{equation}
\end{enumerate}
Comparing Equations (27)-(29) with Equations (13)-(15), we obtain

{\bf Corollary 10} In $Sh[G]$, if we replace $w$ by $1-d^2$, and
replace $a$ by $(Y/B)^{n_a}$ for every chain $a$, where $n_a$ is
the width of the sheaf $a$, then we have
\begin{equation}
Q[G_s]={B^{\sum_{a\in E(G)}n_a}\over d^{p-2k+1}}Sh[G],
\end{equation}
where $p$ is the number of vertices of $G$, and $k$ is the number
of the connected components of $G$.

{\bf Note2} When $t=z_1=z_2=d$,
$$
\aligned W(G)(d,d,d)&=d^{k(G)-1}\sum_{S\subset E(G)}\{\prod_{e\in
S}x_{c(e)}\}\{\prod_{e\notin S}y_{c(e)}\}d^{k<S>-k(G)}d^{n<S>}\\
&=d^{k(G)-1}\sum_{S\subset E(G)}\{\prod_{e\in
S}x_{c(e)}\}\{\prod_{e\notin S}y_{c(e)}\}d^{k<S>-k(G)}d^{(|S|-|V(G)|+k<S>)}\\
&=d^{(-k(G)-|V(G)|)}\sum_{S\subset E(G)}\{\prod_{e\in
S}x_{c(e)}\}\{\prod_{e\notin S}y_{c(e)}\}d^{(|S|+2k<S>)}.
\endaligned
$$

\section{\bf\large Kauffman brackets of rational links}
\noindent

In this section we use the relation built in the last section and
transfer matrix approach to provide an explicit formula for the
Kauffman brackets of general rational links $m_1m_2\cdots m_k$,
where $m_1,m_2$,$\cdots$,$m_k$ are all non-zero integers. There
are two kinds of rational links according to the parity of $k$.
When $k$ is odd we call $m_1m_2\cdots m_k$ horizontal rational
link, and when $k$ is even we call $m_1m_2\cdots m_k$ vertical
rational link. Horizontal rational link $m_1m_2\cdots m_6$ and
vertical rational link $m_1m_2\cdots m_7$ are shown in Figure 4
(upper) respectively. The horizontal and vertical boxes containing
integers $m_i$ represent $m_i$ half-twists shown in Figure 4
(lower).

\begin{figure}[htbp]
\centering
\includegraphics[height=9cm]{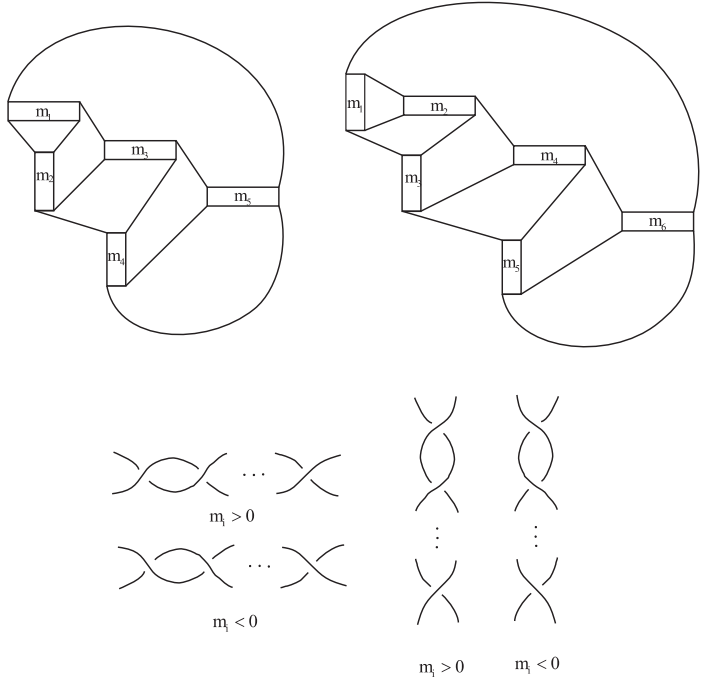}
\caption{Two kinds of rational links}
\end{figure}

Horizontal and vertical rational link diagrams both end with a
horizontal box, and begin with a horizontal box and a vertical box
respectively. The graphs associated with the two kinds of rational
links are shown in Figure 5. We denote them by $HF_n$ and $VF_n$
respectively.

\begin{figure}[htbp]
\centering
\includegraphics[height=6cm]{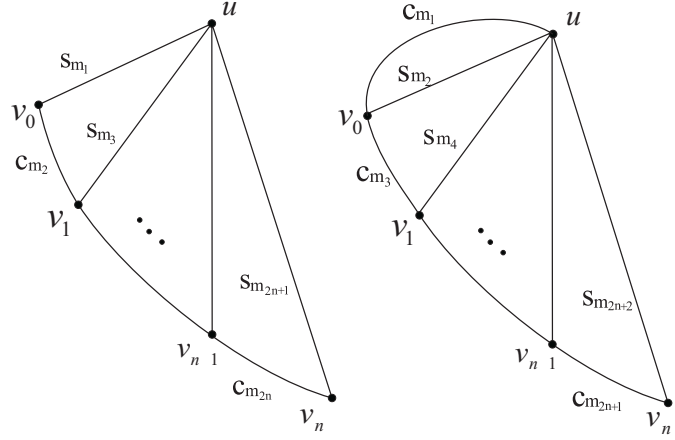}
\caption{Graphs $HF_n$ and $VF_n$ corresponding to rational links}
\end{figure}

{\bf Theorem 11} Let $m_1m_2\cdots m_{2n}m_{2n+1}$ ($n\geq 1$) be
a horizontal rational link.  Let $B=A^{-1}$, $d=-A^2-A^{-2}$,
$X=-A^{-3}$ and $Y=-A^3$. Then
$$
<m_1m_2\cdots m_{2n}m_{2n+1}>=d^{-n-1}X_0^T\{\prod_{i=1}^nA_i\}J,
$$
where $X_0^T=(B^{m_1}d^2,Y^{m_1}-B^{m_1})$,
$$A_i=\left(
\begin{array}{cc}
X^{m_{2i}}B^{m_{2i+1}}d&X^{m_{2i}}(Y^{m_{2i+1}}-B^{m_{2i+1}})d^{-1}\\
(X^{m_{2i}}-A^{m_{2i}})B^{m_{2i+1}}d&(X^{m_{2i}}-A^{m_{2i}})(Y^{m_{2i+1}}-B^{m_{2i+1}})d^{-1}+A^{m_{2i}}Y^{m_{2i+1}}d
\end{array}
\right)
$$
and $J=(1,1)^T$.

{\bf Proof:} For simplicity, we use $W(G)$ to denote
$\sum_{S\subset E(G)}\{\prod_{e\in S}x_{c(e)}\}\{\prod_{e\notin
S}y_{c(e)}\}d^{(|S|+2k<S>)}$.

We denote by $\Lambda_n$ the graph induced by the two edges
$c_{m_{2n}}$ and $s_{m_{2n+1}}$, and view $HF_n$ as the union of
$HF_{n-1}$ and $\Lambda_n$. Note that each $S\subset E(HF_n)$ can
be written as $S=B\cup C$ with $B\subset E(HF_{n-1})$ and
$C\subset E(\Lambda_{n})$. We have
\begin{equation}
\aligned W(HF_n)&=\sum_{S\subset E(HF_n)}\{\prod_{e\in
S}x_{c(e)}\}\{\prod_{e\notin
S}y_{c(e)}\}d^{|S|+2k<S>}\\
&=\sum_{B\cup C}\{\prod_{e\in B}x_{c(e)}\}\{\prod_{e\in
C}x_{c(e)}\}\{\prod_{e\notin B}y_{c(e)}\}\{\prod_{e\notin
C}y_{c(e)}\}d^{|B|+|C|+2(k<B>+\delta(B,C))}\\
&=\sum_{B}\{\prod_{e\in B}x_{c(e)}\}\{\prod_{e\notin
B}y_{c(e)}\}d^{|B|+2k<B>}\sum_{C}\{\prod_{e\in
C}x_{c(e)}\}\{\prod_{e\notin C}y_{c(e)}\}d^{|C|+2\delta(B,C)},
\endaligned
\end{equation}
where $B\subset E(HF_{n-1})$, $C\subset E(\Lambda_{n})$, $e\notin
B$ means $e\in E(HF_{n-1})-B$, $e\notin C$ means $e\in
E(\Lambda_{n})-C$, $k<B>$ is the number of connected components of
$<B>$, the spanning subgraph with edge set $B$ of $HF_{n-1}$, and
$\delta(B,C)=k<B\cup C>-k<B>$, where $k<B\cup C>$ is the number of
connected components of $<B\cup C>$, the spanning subgraph with
edge set $B\cup C$ of $HF_{n}$.

Note that $S\subset E(HF_n)$ and $B\subset E(HF_{n-1})$, we denote
by $s_1$ the state of $S$ (resp. $B$) if $u$ and $v_n$ (resp.
$v_{n-1}$) are disconnected in $<S>$ (resp. $<B>$), and by $s_2$
the state of $S$ (resp. $B$) if $u$ and $v_n$ (resp. $v_{n-1}$)
are connected in $<S>$ (resp. $<B>$). It is not difficult to see
that $\delta(B,C)$ can be computed from the knowledge of the state
of $B$ and $C$. We list them in Table 1. The last two columns of
the table are the state of $S=B\cup C$ and the contribution to the
second sum of Equation (31).

\begin{table}[htp]
\centering
\begin{tabular}{|c|c|c|c|c|}\hline
   the state of $B$  & $C$ &  $\delta(B,C)$& the state of $S$& contribution\\ \hline
   $s_1$ &   $\emptyset$                           &1  &$s_1$ &$y_{c_{m_{2n}}}y_{s_{m_{2n+1}}}d^2$\\ \hline
   $s_1$ &   $\{s_{m_{2n+1}}\}$                    &0  &$s_2$ &$y_{c_{m_{2n}}}x_{s_{m_{2n+1}}}d$\\ \hline
   $s_1$ &   $\{c_{m_{2n}}\}$                      &0  &$s_1$ &$x_{c_{m_{2n}}}y_{s_{m_{2n+1}}}d$\\ \hline
   $s_1$ &   $\{s_{m_{2n+1}},c_{m_{2n}}\}$         &-1 &$s_2$ &$x_{c_{m_{2n}}}x_{s_{m_{2n+1}}}$  \\  \hline
   $s_2$ &   $\emptyset$                           &1  &$s_1$ &$y_{c_{m_{2n}}}y_{s_{m_{2n+1}}}d^2$\\ \hline
   $s_2$ &   $\{s_{m_{2n+1}}\}$                    &0  &$s_2$ &$y_{c_{m_{2n}}}x_{s_{m_{2n+1}}}d$        \\ \hline
   $s_2$ &   $\{c_{m_{2n}}\}$                      &0  &$s_2$ &$x_{c_{m_{2n}}}y_{s_{m_{2n+1}}}d$        \\ \hline
   $s_2$ &   $\{s_{m_{2n+1}},c_{m_{2n}}\}$         &0  &$s_2$ &$x_{c_{m_{2n}}}x_{s_{m_{2n+1}}}d^2$       \\ \hline
\end{tabular}
\caption{The states of $B$ and $S$ and the contributions}
\end{table}

Let $A_n$ be $2\times 2$ matrix with $(i,j)$-entry $a_{i,j}^n$,
where $a_{i,j}^n$ is the sum of contributions from initial state
(the state of $B$) $s_i$ to final state $s_j$ (the state of $S$),
that is,
$$
\aligned
&a_{1,1}^n=y_{c_{m_{2n}}}y_{s_{m_{2n+1}}}d^2+x_{c_{m_{2n}}}y_{s_{m_{2n+1}}}d,\\
&a_{1,2}^n=y_{c_{m_{2n}}}x_{s_{m_{2n+1}}}d+x_{c_{m_{2n}}}x_{s_{m_{2n+1}}},\\
&a_{2,1}^n=y_{c_{m_{2n}}}y_{s_{m_{2n+1}}}d^2{\rm \ \ and}\\
&a_{2,2}^n=x_{c_{m_{2n}}}y_{s_{m_{2n+1}}}d+y_{c_{m_{2n}}}x_{s_{m_{2n+1}}}d+x_{c_{m_{2n}}}x_{s_{m_{2n+1}}}d^2.\\
\endaligned
$$

We divide $B's$ into two subclasses: one is the $B's$ with state
$s_1$ denoted by $B:s_1$, and the other is the $B's$ with state
$s_2$ denoted by $B:s_2$. Then we have

$$
\aligned W(HF_n)&=\sum_{B}\{\prod_{e\in
B}x_{c(e)}\}\{\prod_{e\notin
B}y_{c(e)}\}d^{|B|+2k<B>}\sum_{C}\{\prod_{e\in
C}x_{c(e)}\}\{\prod_{e\notin C}y_{c(e)}\}d^{|C|+2\delta(B,C)}\\
&=\sum_{B:s_1}\{\prod_{e\in B}x_{c(e)}\}\{\prod_{e\notin
B}y_{c(e)}\}d^{|B|+2k<B>}\sum_{C}\{\prod_{e\in
C}x_{c(e)}\}\{\prod_{e\notin C}y_{c(e)}\}d^{|C|+2\delta(B,C)}\\
&\ \ \ \ +\sum_{B:s_2}\{\prod_{e\in B}x_{c(e)}\}\{\prod_{e\notin
B}y_{c(e)}\}d^{|B|+2k<B>}\sum_{C}\{\prod_{e\in
C}x_{c(e)}\}\{\prod_{e\notin C}y_{c(e)}\}d^{|C|+2\delta(B,C)}\\
&=\sum_{B:s_1}\{\prod_{e\in B}x_{c(e)}\}\{\prod_{e\notin
B}y_{c(e)}\}d^{|B|+2k<B>}(a_{1,1}^n+a_{1,2}^n)\\
&\ \ \ \ +\sum_{B:s_2}\{\prod_{e\in B}x_{c(e)}\}\{\prod_{e\notin
B}y_{c(e)}\}d^{|B|+2k<B>}(a_{2,1}^n+a_{2,2}^n).\\
\endaligned
$$

Note that the above two summands.  We denoted them by $s_1^{n-1}$
and $s_2^{n-1}$ respectively, that is,
$$
\aligned
&s_1^{n-1}=\sum_{B:s_1}\{\prod_{e\in
B}x_{c(e)}\}\{\prod_{e\notin
B}y_{c(e)}\}d^{|B|+2k<B>}\\
&s_2^{n-1}=\sum_{B:s_2}\{\prod_{e\in B}x_{c(e)}\}\{\prod_{e\notin
B}y_{c(e)}\}d^{|B|+2k<B>}
\endaligned
$$

On the other hand,
$$
\aligned W(HF_n)&=\sum_{S\subset E(HF_n)}\{\prod_{e\in
S}x_{c(e)}\}\{\prod_{e\notin S}y_{c(e)}\}d^{|S|+2k<S>}\\
&=\sum_{S:s_1}\{\prod_{e\in S}x_{c(e)}\}\{\prod_{e\notin
S}y_{c(e)}\}d^{|S|+2k<S>}+\sum_{S:s_2}\{\prod_{e\in
S}x_{c(e)}\}\{\prod_{e\notin S}y_{c(e)}\}d^{|S|+2k<S>}.\\
\endaligned
$$

If we denote by $s_1^{n}$ and $s_2^{n}$ the above two summands,
then $s_1^{n}=s_1^{n-1}a_{1,1}^n+s_2^{n-1}a_{2,1}^n$ and
$s_2^{n}=s_1^{n-1}a_{1,2}^n+s_2^{n-1}a_{2,2}^n$. Thus we have
$$
\aligned
W(HF_n)&=s_1^{n}+s_2^{n}\\
&=(s_1^{n},s_2^{n})J\\
&=(s_1^{n-1},s_2^{n-1})A_nJ\\
&=(s_1^{n-2},s_2^{n-2})A_{n-1}A_nJ\\
&=\cdots\\
&=(s_1^{0},s_2^{0})\{\prod_{i=1}^nA_i\}J,
\endaligned
$$
where
$$
\aligned
&s_1^{0}=\sum_{S\subset E(HF_0):s_1}\{\prod_{e\in
S}x_{c(e)}\}\{\prod_{e\notin
S}y_{c(e)}\}d^{|S|+2k<S>},\\
&s_2^{0}=\sum_{S\subset E(HF_0):s_2}\{\prod_{e\in
S}x_{c(e)}\}\{\prod_{e\notin S}y_{c(e)}\}d^{|S|+2k<S>},
\endaligned
$$

$A_i= \left(
\begin{array}{cc}
a_{1,1}^i&a_{1,2}^i\\
a_{2,1}^i&a_{2,2}^i
\end{array}
\right)$ with
$$
\aligned
&a_{1,1}^i=y_{c_{m_{2i}}}y_{s_{m_{2i+1}}}d^2+x_{c_{m_{2i}}}y_{s_{m_{2i+1}}}d,\\
&a_{1,2}^i=y_{c_{m_{2i}}}x_{s_{m_{2i+1}}}d+x_{c_{m_{2i}}}x_{s_{m_{2i+1}}},\\
&a_{2,1}^i=y_{c_{m_{2i}}}y_{s_{m_{2i+1}}}d^2{\rm \ \ and}\\
&a_{2,2}^i=x_{c_{m_{2i}}}y_{s_{m_{2i+1}}}d+y_{c_{m_{2i}}}x_{s_{m_{2i+1}}}d+x_{c_{m_{2i}}}x_{s_{m_{2i+1}}}d^2,
\endaligned
$$
and $J=(1,1)^T$.

 Note that $HF_0$ is the graph induced by the edge
$s_{m_1}$, the subset of $E(HF_0)$ with state $s_1$ is $\emptyset$
and the subset of $E(HF_0)$ with state $s_2$ is $\{s_{m_1}\}$. We
obtain $(s_1^{0},s_2^{0})=(y_{s_{m_1}}d^4,x_{s_{m_1}}d^3)$.

By Theorem 8, let $x_{c_n}=A^n, y_{c_n}={X^n-A^n\over d}
x_{s_n}={Y^n-B^n\over d}$ and $y_{s_n}=B^n$, we obtain
$(s_1^{0},s_2^{0})=(B^{m_1}d^4,(Y^{m_1}-B^{m_1})d^2)=d^2(B^{m_1}d^2,Y^{m_1}-B^{m_1})$
and
$$
\aligned
&a_{1,1}^i=X^{m_{2i}}B^{m_{2i+1}}d,\\
&a_{1,2}^i=X^{m_{2i}}(Y^{m_{2i+1}}-B^{m_{2i+1}})d^{-1},\\
&a_{2,1}^i=(X^{m_{2i}}-A^{m_{2i}})B^{m_{2i+1}}d{\rm \ \ and}\\
&a_{2,2}^i=(X^{m_{2i}}-A^{m_{2i}})(Y^{m_{2i+1}}-B^{m_{2i+1}})d^{-1}+A^{m_{2i}}Y^{m_{2i+1}}d.\\
\endaligned
$$
The theorem is proved $\Box$

{\bf Theorem 12} Let $m_1m_2\cdots m_{2n+1}m_{2n+2}$ ($n\geq 1$)
be a vertical rational link. Let $B=A^{-1}$, $d=-A^2-A^{-2}$,
$X=-A^{-3}$ and $Y=-A^3$. Then
$$
<m_1m_2\cdots
m_{2n+1}m_{2n+2}>=d^{-n-2}X_0^T\{\prod_{i=1}^nA_i\}J,
$$
where
$X_0^T=((X^{m_1}-A^{m_1})B^{m_2}d^2,A^{m_1}Y^{m_2}d^2+(X^{m_1}-A^{m_1})(Y^{m_2}-B^{m_2}))$,
$$A_i=\left(
\begin{array}{cc}
X^{m_{2i+1}}B^{m_{2i+2}}d&X^{m_{2i+1}}(Y^{m_{2i+2}}-B^{m_{2i+2}})d^{-1}\\
(X^{m_{2i+1}}-A^{m_{2i+2}})B^{m_{2i+2}}d&(X^{m_{2i+1}}-A^{m_{2i+1}})(Y^{m_{2i+2}}-B^{m_{2i+2}})d^{-1}+A^{m_{2i+1}}Y^{m_{2i+2}}d
\end{array}
\right)
$$
and $J=(1,1)^T$.

{\bf Proof:} The proof is similar to that of Theorem 11. The only
difference is that $VF_0$ is the graph induced by two edges
$c_{m_1}$ and $s_{m_2}$, the subset of $E(VF_0)$ with state $s_1$
is $\emptyset$ and the subset of $E(HF_0)$ with state $s_2$ is
$\{c_{m_1}\}$,$\{s_{m_2}\}$ and $\{c_{m_1},s_{m_2}\}$. Thus
$$
\aligned &s_1^0=y_{c_{m_1}}y_{s_{m_2}}d^4,\\
&s_2^0=x_{c_{m_1}}y_{s_{m_2}}d^3+y_{c_{m_1}}x_{s_{m_2}}d^3+x_{c_{m_1}}x_{s_{m_2}}d^4.
\endaligned
$$
After replacement, we have
$$
\aligned &s_1^0=(X^{m_1}-A^{m_1})B^{m_2}d^3,\\
&s_2^0=A^{m_1}Y^{m_2}d^3+(X^{m_1}-A^{m_1})(Y^{m_2}-B^{m_2})d.
\endaligned
$$
$\Box$

There are only $m_1m_2$ (i.e. generalized twist link) shown in
Figure 3 and $m_1$ (i.e. ($m_1$,2)-torus link) shown in Figure 6
left for us to compute.
\begin{figure}[htbp]
\centering
\includegraphics[height=2.5cm]{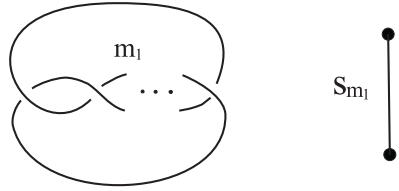}
\caption{($m_1$,2)-torus link and its associated graph}
\end{figure}

\begin{enumerate}
\item Kauffman bracket of rational link $m_1$.

The graph associated with $m_1$ is the graph induced by an edge
$e_1$, which represents a sheaf with width $m_1$ and its
$W$-polynomial is $x_{c(e_1)}+z_1y_{c(e_1)}.$ By Theorem 8, the
Kauffman bracket of $m_1$ is
$$
{Y^{m_1}-B^{m_1}\over d}+B^{m_1}d={1\over
-A^2-A^{-2}}\{(-A^3)^{m_1}-A^{-m_1}\}+A^{-m_1}(-A^2-A^{-2}).
$$
\item Kauffman bracket of rational link $m_1m_2$.

The graph associated with $m_1m_2$ is the graph with two vertices
by two parallel edges, say, $e_1$ and $e_2$ (see Figure 3, right).
By the definition or Theorem 4, the $W$-polynomial of this graph
is
$$
(x_{c(e_1)}+z_1y_{c(e_1)})y_{c(e_2)}+(x_{c(e_1)}z_2+y_{c(e_1)})x_{c(e_2)}.
$$

Note that the length of the chain $e_1$ is $m_1$ and the width of
the sheaf $e_2$ is $m_2$, by Theorem 8, the Kauffman bracket of
$m_1m_2$ is
$$
(A^{m_1}+d{X^{m_1}-A^{m_1}\over
d})B^{m_2}+(A^{m_1}d+{X^{m_1}-A^{m_1}\over
d}){Y^{m_2}-B^{m_2}\over d}.
$$
After reduction, we obtain
\begin{equation}
\aligned
<m_1m_2>=&A^{m_1-m_2}\left\{(-A^{-4})^{m_1}+(-A^4)^{m_2}-1\right\}+\\
&{A^{m_1-m_2}\over
A^{-4}+2+A^4}\left\{(-A^{-4})^{m_1-m_2}-(-A^{-4})^{m_1}-(-A^{4})^{m_2}+1\right\}.
\endaligned
\end{equation}
\end{enumerate}
\section{\bf\large Further examples}
\noindent

It is clear we can also use Theorem 8 to compute the Kauffman
brackets of other link families of rational links type. For
example, using Corollary 9 the explicit formulae of Kauffman
brackets of pretzel links are obtained by the present author in
[8]. In this section we provide another example. We denote by
$L(m_1,m_2,m_3)$ the link shown in Figure 7 when $m_i>0$ for
$i=1,2,3$. When $m_i<0$, the link diagram is obtained by replacing
the corresponding half-twists by its mirror image.

\begin{figure}[htbp]
\centering
\includegraphics[height=4cm]{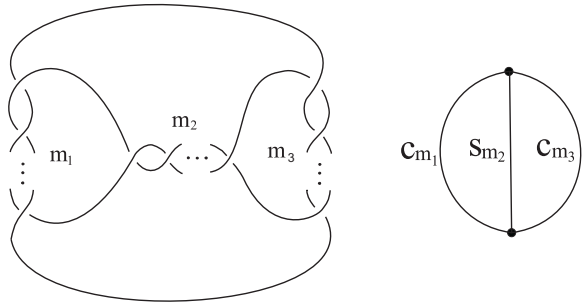}
\caption{The link $L(m_1,m_2,m_3)$ and its associated graph
$\Theta$}
\end{figure}

The graph $\Theta$ associated with $L(m_1,m_2,m_3)$ has three
spanning trees with edge sets $\{c_{m_1}\},\{s_{m_2}\}$ and
$\{c_{m_3}\}$ respectively. By Theorem 4, we have
$$
\aligned
W(\Theta)(d,d,d)&=\sum_F\{\prod_{IA}(x_{c(e)}+dy_{c(e)})\}\{\prod_{EA}(x_{c(e)}d+y_{c(e)})\}\{\prod_{II}x_{c(e)}\}\{\prod_{EI}y_{c(e)}\}\\
&=(x_{c_{m_1}}+dy_{c_{m_1}})y_{s_{m_2}}y_{c_{m_3}}+(x_{c_{m_1}}d+y_{c_{m_1}})x_{s_{m_2}}y_{c_{m_3}}\\
&\ \ \ \
+(x_{c_{m_1}}d+y_{c_{m_1}})(x_{s_{m_2}}d+y_{s_{m_2}})x_{c_{m_3}}.
\endaligned
$$
By Theorem 8, we obtain the Kauffman brackets of $L(m_1,m_2,m_3)$
$$
\aligned
<L(m_1,m_2,m_3)>&=X^{m_1}B^{m_2}(X^{m_3}-A^{m_3})/d+\\
&(A^{m_1}d+(X^{m_1}-A^{m_1})/d)(Y^{m_2}-B^{m_2})/d(X^{m_3}-A^{m_3})/d+\\
&((d^2-1)A^{m_1}+X^{m_1})Y^{m_2}A^{m_3}.
\endaligned
$$

 \vskip2cm
\begin{center}
{\bf\large References}
\end{center}
\vskip0.4cm
\begin{description}
\item{[1]}\ L.H. Kauffman, State models and the Jones polynomial,
Topology 26 (1987) 395-407. \item{[2]}\ V.F.R. Jones, A polynomial
invariant for knots via Von Neumann algebras. Bull. Am. Math. Soc.
12 (1985) 103-112.
\item{[3]}\ R.A. Landvoy, The Jones polynomial
of pretzel knots and links, Topology and its Applications 83
(1998) 135-147.
\item{[4]}\ P.M.G. Manch$\acute{o}$n, On the
Kauffman bracket of pretzel links, Marie Curie Fellowships Annals,
Second Volume, 2003, 118-122.

\item{[5]}\ L.H. Kauffman, New invariants in the theory of knots,
Am. Math. Monthly 95 (1988) 195-242.
\item{[6]}\ L.H. Kauffman, A
Tutte polynomial for signed graphs, Discrete Applied Mathematics
(1989) 105-127.

\item{[7]}\ X.A. Jin, F.J.
Zhang, Zeros of the Jones polynomials for families of pretzel,
Physica A 328 (2003) 391-408.
\item{[8]}\ X.A. Jin, F.J. Zhang, The
Kauffman brackets for equivalence classes of links, Advances in
applied mathematics 34 (2005) 47-64.

\item{[9]}\ B. Bollob$\acute{a}$s, A Tutte polynomial for colored
graphs, Combinatorics, Probability and Computing 8 (1999) 45-93.

\item{[10]}\ J.H. Conway, An enumeration of knots and links, and
some of their algebraic properties, in Computational Problems in
Abstract Alegebra, Pergamon Press, New York (1969) 329-358.

\item{[11]}\ J.R. Goldman, L.H. Kauffman, Rational tangles,
Advance in applied mathematics 18 (1997) 300-332.

 \item{[12]}\ B. Bollob$\acute{a}$s, Modern Graph
Theory, Chapter 10, Springer, 2001.

\item{[13]}\ W.T. Tutte, A contribution to the theory of chromatic
polynomials, Can. J. Math.6 (1954) 80-91.

\item{[14]}\ R.C. Read, E.G. Whitehead Jr., Chromatic polynomials
of homeomorphism classes of graphs, Discrete Math. 204 (1999)
337-356.

\item{[15]}\ L. Traldi, A dichromatic polynomial for weighted
graphs and link diagrams, Proceedings of the American Mathematical
Society, Volume 106, Number 1, May 1989, 279-286.

\item{[16]}\ L. Traldi, Chain polynomials and Tutte polynomials,
Discrete Math. 248 (2002) 279-282.

\end{description}

\end{document}